\documentclass{article}
\usepackage{geometry,enumerate}
\usepackage{amsmath, amssymb,amsthm,mathtools}
\usepackage{xypic}
\usepackage[all]{xy}
\usepackage{graphics}
\usepackage{comment}
\usepackage{cleveref}
\usepackage{multicol}\columnseprule=0.4pt
\usepackage{color}
\title{Connectedness of quasi-hereditary structures}
\author{Yuichiro Goto}
\date{}
\newtheorem{df}{Definition}[section]
\newtheorem{eg}[df]{Example}
\newtheorem{thm}[df]{Theorem}
\newtheorem{prop}[df]{Proposition}
\newtheorem{lem}[df]{Lemma}
\newtheorem{cor}[df]{Corollary}
\newtheorem{rmk}[df]{Remark}

\newenvironment{thmbis}[1]
{\renewcommand{\thedf}{\ref{#1}$'$}%
\addtocounter{df}{-1}%
\begin{thm}}
{\end{thm}}
\newenvironment{lembis}[1]
{\renewcommand{\thedf}{\ref{#1}$'$}%
\addtocounter{df}{-1}%
\begin{lem}}
{\end{lem}}

\crefname{df}{Definition}{Definitions}
\crefname{eg}{Example}{Examples}
\crefname{thm}{Theorem}{Theorems}
\crefname{prop}{Proposition}{Propositions}
\crefname{lem}{Lemma}{Lemmas}
\crefname{cor}{Corollary}{Corollaries}
\crefname{rmk}{Remark}{Remarks}

\newcommand{\dra}{{\dashrightarrow}}

\newcommand{\bm}[1]{{\mbox{\boldmath $#1$}}}
\newenvironment{smatrix}{ \begin{smallmatrix}} {\end{smallmatrix}}

\DeclareMathOperator{\rad}{rad}

\DeclareMathOperator{\kernel}{Ker}
\DeclareMathOperator{\image}{Im}
\DeclareMathOperator{\Coker}{Coker}

\DeclareMathOperator{\module}{mod}

\DeclareMathOperator{\Hom}{Hom}

\DeclareMathOperator{\Ext}{Ext}
\DeclareMathOperator{\tr}{tr}

\begin{document}

\maketitle \vspace{3mm}

\begin{abstract}
Dlab and Ringel showed that algebras being quasi-hereditary in all orders for indices of primitive idempotents becomes hereditary. 
So, we are interested in for which orders a given quasi-hereditary algebra is again quasi-hereditary. 
As a matter of fact, we consider permutations of indices, and if the algebra with permuted indices is quasi-hereditary, then we say that this permutation gives a quasi-hereditary structure.

In this article, we first give a criterion for adjacent transpositions giving quasi-hereditary structures, in terms of homological conditions of standard or costandard modules over a given quasi-hereditary algebra. 
Next, we consider those which we call connectedness of quasi-hereditary structures. 
The definition of connectedness can be found in \Cref{def of connected}.
We then show that any two quasi-hereditary structures are connected, which is our main result. 
By this result, once we know that there are two quasi-hereditary structures, then permutations in some sense lying between them give also quasi-hereditary structures.
\end{abstract} 

\section{Introduction}
Quasi-hereditary algebras, introduced by Cline, Parshall and Scott, generalize hereditary algebras.
Moreover Dlab and Ringel showed in Theorem 1 of \cite{DR1} that if an algebra is quasi-hereditary in all orders, it becomes hereditary, and vice versa.
From this point of view, we study quasi-hereditary structures for a given algebra.
Recently, there are two results on quasi-hereditary structures.
Coulembier showed in \cite{C} that a quasi-hereditary algebra with simple preserving duality has only one quasi-hereditary structure.
Flores, Kimura and Rognerud gave a method of counting the number of quasi-hereditary structures for a path algebras of Dynkin types in \cite{FKR}.
In their papers, the quasi-hereditary structure was defined by an equivalent class of partial orders with some relations.
However in this article, we define it by using a total order without using equivalent classes.
Thus, our results are in the nature different from them and can not be derived from their results.
Moreover we will use permutations instead of total orders when considering quasi-hereditary structures.

Fix a finite dimensional algebra $A$ with pairwise orthogonal primitive idempotents $e_1,\dots,e_n$.
Let $\mathfrak{S}_n$ be the symmetric group on $n$ letters and $\sigma \in \mathfrak{S}_n$.
We say that $\sigma$ gives a quasi-hereditary structure of $A$ if $e_{\sigma^{-1}(1)},\dots,e_{\sigma^{-1}(n)}$ with this order, called the $\sigma$-order, make $A$ quasi-hereditary, see \Cref{def of qha}.
Let $X$ be a certain property defined for elements of $\mathfrak{S}_n$.
If $\sigma,\tau \in \mathfrak{S}_n$ have the property $X$, and if there is a decomposition 
\[\tau \sigma^{-1} = \sigma_{i_l} \cdots \sigma_{i_1}\]
into the product of adjacent transpositions such that all $\sigma_{i_k} \cdots \sigma_{i_1} \sigma$ also have the property $X$ for $1 \leq k \leq l$, then we say that $\sigma$ and $\tau$ are connected with respect to $X$.
For precise definitions, see \Cref{def of connected}.

This paper is organized as follows.
\Cref{Preliminary} recalls the definition of quasi-hereditary algebras and introduces their properties.
In \Cref{twistability}, we discuss when a quasi-hereditary structure gives another quasi-hereditary structure.
For a quasi-hereditary algebra, by using homological conditions of standard or costandard modules over it, we give a criterion for adjacent transpositions giving quasi-hereditary structures.
\Cref{Connectedness} shows the connectedness of two permutations giving quasi-hereditary structures.
The following is the main result of this article.
\begin{thm}[\Cref{connected}]
    Any two permutations giving quasi-hereditary structures are connected.
\end{thm}
Moreover, we give a method to obtain some permutations giving quasi-hereditary structures which induce the connectedness of the two.
\begin{cor}[\Cref{decomp of connected}]     
    Let $\sigma,\tau$ give quasi-hereditary structures with $\sigma \neq \tau$.\\
    For $k = 1,2,\dots$, inductively take a minimal element $i_{k}$ with respect to the $(\sigma_{i_{k-1}} \cdots \sigma_{i_1} \sigma)$-order satisfying $\sigma_{i_{k-1}} \cdots \sigma_{i_1} \sigma(i_{k}) \neq n$ and
    \[\tau(i_{k}) > \tau\sigma^{-1} \sigma_{i_1} \cdots \sigma_{i_{k-1}}(\sigma_{i_{k-1}} \cdots \sigma_{i_1} \sigma (i_{k})+1).\]
    We take $i_1,i_2,\dots,i_k$ until those elements satisfying the above exist.
    If there is no $i_{k+1}$ satisfying the above, then we do not take $i_{k+1}$ and put $l=k$.
    Then the product $\sigma_{i_l} \cdots \sigma_{i_1}$ is a decomposition of $\tau\sigma^{-1}$ inducing the connectedness of $\sigma$ and $\tau$.
\end{cor}

{\bfseries Acknowledgment}\\
I would like to thank prof. Katsuhiro Uno for many helpful discussions with him and improving this article.
I am also grateful to Takahide Adachi, Aaron Chan, Yuta Kimura and Mayu Tsukamoto for their valuable comments and discussions.

\section{Preliminary}\label{Preliminary}
Throughout this article, let $K$ be an algebraically closed field, $A$ a finite dimensional $K$-algebra with pairwise orthogonal primitive idempotents $e_1,\dots,e_n$, and let $\Lambda=\{1,\dots,n\}$.
For $i \in \Lambda$, we denote $P(i)=e_i A$ the indecomposable projective module, $S(i)$ the top of $P(i)$, and $I(i)$ the injective envelope of $S(i)$.
The category of finitely generated right $A$-modules is denoted by $\module A$ and call its objects just $A$-modules.
The standard $K$-dual $\Hom_K(-,K)$ is denoted by $D$. 
For an $A$-module $M$, we write the isomorphism class of $M$ by $[M]$ and the Jordan-H\"older multiplicity of $S(i)$ in $M$ by $[M:S(i)]$.
Write $e \in \mathfrak{S}_n$ the trivial permutation and $\sigma_i=(i,i+1) \in \mathfrak{S}_n$ adjacent transpositions for $1 \leq i \leq n-1$.

In this section, we recall the definition and the well-known properties of quasi-hereditary algebras. 

\begin{df}\label{def of qha}
    Let $A$ be an algebra as above and $\sigma \in \mathfrak{S}_n$.
\begin{enumerate}[$(1)$]
    \item For $i,j \in \Lambda$, we write $i <^\sigma j$ as $\sigma(i) < \sigma(j)$.
    Then there is the total order 
    \[\sigma^{-1}(1) <^\sigma \sigma^{-1}(2) <^\sigma \cdots <^\sigma \sigma^{-1}(n)\] 
    over $\Lambda$, and it is called the {\bfseries ${\bm \sigma}$-order}.
    \item For each $i \in \Lambda$, the $A$-module $\Delta^\sigma(i)$, called the {\bfseries standard module} with respect to the $\sigma$-order, is defined by the maximal factor module of $P(i)$ having only composition factors $S(j)$ with $\sigma(j) \leq \sigma(i)$.
    Moreover we will write the set $\{\Delta^\sigma(1), \dots, \Delta^\sigma(n)\}$ by $\Delta^\sigma$. 
    \item Dually, the $A$-module $\nabla^\sigma(i)$, called the {\bfseries costandard module} with respect to the $\sigma$-order, is defined by the maximal submodule of $I(i)$ having only composition factors $S(j)$ with $\sigma(j) \leq \sigma(i)$.
    Denote the set $\{\nabla^\sigma(1), \dots, \nabla^\sigma(n)\}$ by $\nabla^\sigma$. 
    \item We say that an $A$-module $M$ has a {\bfseries $\Delta^\sigma$-filtration} (resp. a {\bfseries $\nabla^\sigma$-filtration}) if there is a sequence of submodules 
    \[M=M_1 \supset M_2 \supset \dots \supset M_{m+1} = 0\] 
    such that for each $1 \leq k \leq m$, $M_k / M_{k+1} \cong \Delta^\sigma(j)$ (resp. $M_k / M_{k+1} \cong \nabla^\sigma(j)$) for some $j\in\Lambda$.
    \item A pair $(A,\sigma)$ is said to be a {\bfseries quasi-hereditary algebra} provided that the following conditions are satisfied: 
    \begin{enumerate}[$(a)$]
        \item $[\Delta^\sigma(i):S(i)]=1$ for all $i \in \Lambda$.
        \item $A_A$ has a $\Delta^\sigma$-filtration.
    \end{enumerate}
    If this is the case, we say that the permutation $\sigma$ gives a {\bfseries quasi-hereditary structure} of $A$.
\end{enumerate}
\end{df}

\begin{rmk}[\cite{DR1}]\label{qh rmk}
    \begin{enumerate}[$(1)$]
        \item Put $\displaystyle\varepsilon^{\sigma}_{i} = \sum_{\sigma(j) \geq \sigma(i)}e_j$ and $J^{\sigma}_i = A\varepsilon^{\sigma}_{i}A$ for each $i\in\Lambda$.
            Now consider a sequence
            \[A = J^{\sigma}_{\sigma^{-1}(1)} \supset J^{\sigma}_{\sigma^{-1}(2)} \supset \cdots \supset J^{\sigma}_{\sigma^{-1}(n)} \supset 0.\]
            Then $(A,\sigma)$ is quasi-hereditary if and only if the following three conditions holds: 
            \begin{enumerate}
                \item $J^{\sigma}_{\sigma^{-1}(i)}/J^{\sigma}_{\sigma^{-1}(i+1)}$ are projective right $A/J^{\sigma}_{\sigma^{-1}(i+1)}$-modules for all $1 \leq i \leq n-1$.
                \item $\varepsilon^{\sigma}_{\sigma^{-1}(i)}\rad({A/J^{\sigma}_{\sigma^{-1}(i+1)}})\varepsilon^{\sigma}_{\sigma^{-1}(i)} = 0$ for all $1 \leq i \leq n-1$.
                \item $J^{\sigma}_{\sigma^{-1}(n)}$ is a projective right $A$-module and $e_{\sigma^{-1}(n)}\rad{A}e_{\sigma^{-1}(n)}=0$.
            \end{enumerate}
            If $(A,\sigma)$ is quasi-hereditary, then the above sequence of ideals is called a {\bfseries heredity chain}.
        \item Let a sequence
        \[A = J^{\sigma}_{\sigma^{-1}(1)} \supset J^{\sigma}_{\sigma^{-1}(2)} \supset \cdots \supset J^{\sigma}_{\sigma^{-1}(n)} \supset 0,\]
        where $J^{\sigma}_i = A\varepsilon^{\sigma}_{i}A$, be a heredity chain.
        Then for any $i \in \Lambda$, the pairs $(\varepsilon^{\sigma}_{i}A\varepsilon^{\sigma}_{i},\sigma')$ and $(A/J^{\sigma}_i,\overline{\sigma})$ are also quasi-hereditary where $\sigma'(j)=\sigma(j)$ for $j \geq i$ and $\overline{\sigma}(k)=\sigma(k)$ for $k < i$.
    \end{enumerate}
\end{rmk}

For simplicity, we abbreviate $\Delta^e,\nabla^e$ as $\Delta,\nabla$ and call just standard, costandard modules, respectively.
By the definitions of standard and costandard modules, we immediately obtain the following properties.
\begin{lem}[\cite{DR2} Lemmas 1.2, 1.3]\label{DR2}    
    We have the following equalities.
    \begin{enumerate}[$(1)$]
        \item For $i > j$, $\Hom_A(\Delta(i),\Delta(j))=\Hom_A(\nabla(j),\nabla(i))=0$.
        \item For $i \geq j$, $\Ext^1_A(\Delta(i),\Delta(j))=\Ext^1_A(\nabla(j),\nabla(i))=0$.
    \end{enumerate}
\end{lem}

Next we introduce some properties of the standard and costandard modules over a quasi-hereditary algebra.
It is well known that $\{[S(1)],\dots,[S(n)]\}$ forms a $\mathbb{Z}$-basis of Grothendieck group $K_0(\module A)$ of $\module A$.
The following theorem shows that so does the standard modules.

\begin{thm}[\cite{D}A.1(7)]\label{Gro basis}
    Assume $[\Delta(i):S(i)]=1$ for all $i \in \Lambda$.
    Then $\{[\Delta(1)],\dots,[\Delta(n)]\}$ forms a $\mathbb{Z}$-basis of $K_0(\module A)$.
    In particular, for a module having a $\Delta$-filtration, the $\Delta$-filtration factors are unique and do not depend on the choice of $\Delta$-filtrations.
\begin{proof}
    The isomorphism class of each standard module $[\Delta(i)]$ is of the forms 
    \[[\Delta(i)]=[S(i)]+\sum_{j < i}[\Delta(i):S(j)][S(j)].\]
    So we have an equality
    \[[S(i)]=[\Delta(i)]-\sum_{j < i}[\Delta(i):S(j)][S(j)].\]
    Thus we can replace $\{[S(1)],\dots,[S(n)]\}$ by $\{[\Delta(1)],[S(2)],\dots,[S(n)]\}$ as a $\mathbb{Z}$-basis of $K_0(\module A)$.
    Inductively, we obtain the desired result.
\end{proof}
\end{thm} 
We denote the $\Delta$-filtration multiplicity of $\Delta(i)$ in $M$ by $[M:\Delta(i)]$. 

As a useful fact, we show some properties which every pair of neighbor standard modules has.
\begin{lem}[\cite{MO}Lemma 2.]\label{HomExt}
    Assume that $(A,e)$ is a quasi-hereditary algebra.
    Then we have the following equalities.
    \begin{enumerate}[$(1)$]
        \item $\dim\Hom_A(\Delta(i),\Delta(i+1))=\dim\Hom_A(P(i),\Delta(i+1))=[\Delta(i+1):S(i)]$.
        \item $\dim\Ext_A^1(\Delta(i),\Delta(i+1))=\dim\Ext_A^1(\Delta(i),S(i+1))=[P(i):\Delta(i+1)]$.
    \end{enumerate}
\end{lem}
We will denote $H_i=\dim\Hom_A(\Delta(i),\Delta(i+1))$ and $E_i=\dim\Ext_A^1(\Delta(i),\Delta(i+1))$.
Now we remark that the dual statements of \Cref{Gro basis,HomExt} also hold as follows.
\begin{thmbis}{Gro basis}
    Assume that $(A,e)$ is a quasi-hereditary algebra.
    Then $\{[\nabla(1)],\dots,[\nabla(n)]\}$ forms a $\mathbb{Z}$-basis of $K_0(\module A)$.
\end{thmbis}

\begin{lembis}{HomExt}\label{HomExt2}   
    Assume that $(A,e)$ is a quasi-hereditary algebra.
    Then we have the following equalities.
    \begin{enumerate}[$(1)$]
        \item $\overline{H_i} :=\dim\Hom_A(\nabla(i+1),\nabla(i))=\dim\Hom_A(\nabla(i+1),I(i))=[\nabla(i+1):S(i)]$.
        \item $\overline{E_i} :=\dim\Ext_A^1(\nabla(i+1),\nabla(i))=\dim\Ext_A^1(S(i+1),\nabla(i))=[I(i):\nabla(i+1)]$.
    \end{enumerate}
\end{lembis}

For convenience, we introduce the following.
\begin{df}
    We denote by $\tr^{\sigma}_{i}(M)$ the {\bfseries trace} of $\{P(j)\}_{\sigma(j) \geq \sigma(i)}$ in $M$, that is,
    \[\tr^{\sigma}_{i}(M) = \sum\{\image\phi \mid \phi:P(j) \to M,\sigma(j) \geq \sigma(i)\}.\]
\end{df}
Then the standard modules with respect to the $\sigma$-order are written as follows.
\begin{lem}[\cite{DR2} Lemma 1.1']\label{stand trace}
    For any $\sigma \in \mathfrak{S}_n$ and $i \in \Lambda$, we have
    \[\Delta^{\sigma}(i) \cong \frac{P(i)}{\tr^{\sigma}_{i_{\sigma+}}\left(P(i)\right)},\]
    where $\sigma(i_{\sigma+}) = \sigma(i)+1$
\end{lem}
Here, we confirm some properties of the map $\tr^{\sigma}_{i}$.
\begin{lem}\label{trace}
    \begin{enumerate}[$(1)$]
        \item For a homomorphism $f:M \to N$, the map $\tr^{\sigma}_{i}(f):\tr^{\sigma}_{i}(M) \to \tr^{\sigma}_{i}(N)$ makes the following diagram commutative.
        \[\xymatrix{
        \tr^{\sigma}_{i}(M) \ar@{^{(}->}[d] \ar[r]^-{\tr^{\sigma}_{i}(f)} & \tr^{\sigma}_{i}(N) \ar@{^{(}->}[d]\\
        M \ar[r]^-f & N
        }\]
        \item The map $\tr^{\sigma}_{i}$ preserves the surjectivity.
        \item If $M'$ is a submodule of $M$, there is an isomorphism $\tr^{\sigma}_{i}(M/M') \cong \tr^{\sigma}_{i}(M)/(M' \cap \tr^{\sigma}_{i}(M))$.
    \end{enumerate}
    \begin{proof}
        $(1)$ We get $\tr^{\sigma}_{i}(f)$ as the restriction map of $f$ to $\tr^{\sigma}_{i}(M)$.

        $(2)$ By the definitions of $\tr^{\sigma}_{i}$ and projective modules, if $f$ is surjective, then so is $\tr^{\sigma}_{i}(f)$.

        $(3)$ We have the canonical epimorphism $\tr^{\sigma}_{i}(M) \to \tr^{\sigma}_{i}(M/M')$ with kernel $M' \cap \tr^{\sigma}_{i}(M)$, by $(2)$.
    \end{proof}
\end{lem}

\begin{rmk}
    The module $\tr^{\sigma}_{i}A$ is an ideal $J^\sigma_{i}$ given in \Cref{qh rmk}.
    In particular, for each $i\in\Lambda$, the sequence
    \[P(i) = \tr^{\sigma}_{i}P(i) \supset \tr^{\sigma}_{i_{\sigma+}}P(i) \supset \cdots \supset \tr^{\sigma}_{\sigma^{-1}(n)}P(i) \supset 0\]
    satisfies $\tr^{\sigma}_{k}P(i)/\tr^{\sigma}_{k_{\sigma+}}P(i) \cong \Delta^{\sigma}(k)^{[P(i):\Delta^{\sigma}(k)]}$ for all $\sigma(i) \leq \sigma(k) < n$ and $\tr^{\sigma}_{k}P(i) \cong \Delta^{\sigma}(k)^{[P(i):\Delta^{\sigma}(k)]}$ for $\sigma(k) = n$.
    Moreover, this sequence is uniquely defined.
    Here, $k_{\sigma+} = \sigma^{-1}(\sigma(k)+1)$.
\end{rmk}

\section{Twistability}\label{twistability}

Let $(A,\sigma)$ be a quasi-hereditary algebra.
If $(A,\sigma_i\sigma)$ is also quasi-hereditary, then we call the original quasi-hereditary algebra $(A,\sigma)$ to be {\bfseries $\bm i$th-twistable}. 
In this section, we will give the condition on standard or costandard modules equivalent to the $i${\rm th}-twistability for a quasi-hereditary algebra.

First we introduce the homological relations for standard modules and costandard modules. 
\begin{lem}[\cite{MO}Corollary 1]\label{MO}
    Assume that $(A,e)$ is a quasi-hereditary algebra.
    Then the followings hold.
    \begin{enumerate}[$(1)$]
        \item $\Hom_A(\Delta(i),\Delta(i+1)) \cong D\Ext^1_A(\nabla(i+1),\nabla(i))$.
        \item $\Ext^1_A(\Delta(i),\Delta(i+1)) \cong D\Hom_A(\nabla(i+1),\nabla(i))$.
    \end{enumerate}
    In particular, we have $H_i=\overline{E_i}$ and $E_i=\overline{H_i}$.
\end{lem}

\begin{lem}[\cite{AHLU} Theorem 1.1]\label{AHLU}
    Assume that we have $[\Delta(i):S(i)]=1$ for all $i \in \Lambda$.
    Then $A_A$ has a $\Delta$-filtration if and only if $D(_AA)$ has a $\nabla$-filtration.
\end{lem}

\begin{lem}\label{ith-twistable lem}
    Let $(A,e)$ be quasi-hereditary algebra.
    Then the followings are equivalent.
    \begin{enumerate}
        \item $[\Delta^{\sigma_i}(k):S(k)]=1$ for all $k \in \Lambda$.
        \item $E_i\overline{E_i}=0$.
    \end{enumerate}
    \begin{proof}        
        For $k \neq i,i+1$, obviously we have $[\Delta^{\sigma_i}(k):S(k)]=1$ as $\Delta^{\sigma_i}(k)=\Delta(k)$.
        Moreover, the module $\Delta^{\sigma_i}(i+1)$ always satisfies $[\Delta^{\sigma_i}(i+1):S(i+1)]=1$ because it is a factor module of $\Delta(i+1)$.
        So we will claim that $[\Delta^{\sigma_i}(i):S(i)]=1$ if and only if we have $E_i\overline{E_i}=0$.
        The module $\Delta^{\sigma_i}(i)$ is isomorphic to the factor module $P(i)/\tr^{\sigma_i}_{i+2}P(i) = P(i)/\tr^e_{i+2}P(i)$, and hence it is given by the short exact sequence 
        \[0 \to \Delta(i+1)^{E_i} \to \Delta^{\sigma_i}(i) \to \Delta(i) \to 0.\]
        So $[\Delta^{\sigma_i}(i):S(i)] = 1 + E_iH_i$. 
        By using \Cref{MO}, $(A,e)$ is $i${\rm th}-twistable if and only if we have $1 = 1 + E_i\overline{E_i}$, i.e., $E_i\overline{E_i}=0$.
    \end{proof}
\end{lem}

Now we give the condition equivalent to the $i${\rm th}-twistability.
\begin{thm}\label{ith-twistable}
    A quasi-hereditary algebra $(A,e)$ is $i${\rm th}-twistable if and only if one of the following conditions holds:
    \begin{description}
        \item[$(\mathcal{E}_i)$] $E_i=0$ and $\Delta(i+1)$ has a submodule isomorphic to $\Delta(i)^{H_i}$. 
        \item[$(\overline{\mathcal{E}_i})$] $\overline{E_i}=0$ and $\nabla(i+1)$ has a factor module isomorphic to $\nabla(i)^{\overline{H_i}}$. 
    \end{description}
    In particular, if a quasi-hereditary algebra $(A,e)$ satisfies $E_i=\overline{E_i}=0$, then $(A,\sigma_i)$ is also quasi-hereditary with $\Delta^{\sigma_i}=\Delta$ and $\nabla^{\sigma_i}=\nabla$. 
    \begin{proof}
        In the previous lemma, we showed that $[\Delta^{\sigma_i}(k):S(k)]=1$ for all $k \in \Lambda$ if and only if $E_i\overline{E_i}=0$.
        So we concentrate on the situation in which $_AA$ has a $\Delta^{\sigma_i}$-filtration.

        Now show the ``only if'' part.
        We may assume that $E_i=0$ or $\overline{E_i}=0$ from the above discussion.
        If $E_i=0$ and $\overline{E_i}=0$ hold, we have $\Delta^{\sigma_i}=\Delta$ and $\nabla^{\sigma_i}=\nabla$ by \Cref{HomExt,HomExt2}.
        So we first assume that $E_i=0$ and $\overline{E_i} \neq 0$.
        Then by the above arguments, we have $\Delta^{\sigma_i}(k)=\Delta(k)$ for $k\neq i+1$.
        Since $P(i+1)$ has a $\Delta^{\sigma_i}$-filtration by our assumption, so does $\Delta(i+1)$.
        Now there is a factor module of $\Delta(i+1)$ isomorphic to $\Delta^{\sigma_i}(i+1)$ and the kernel of an epimorphism $\Delta(i+1) \to \Delta^{\sigma_i}(i+1)$ is filtered by $\{\Delta^{\sigma_i}(k) \mid \sigma(k) > \sigma(i+1)\}$ by \Cref{DR2}.
        Moreover among the labels of the composition factors of $\Delta(i+1)$, the first and second maximum elements with respect to the $\sigma_i$-order are $i$ and $i+1$, respectively.
        Hence the isomorphism class $[\Delta(i+1)]$ is written by
        \[[\Delta(i+1)] = [S(i+1)] + H_i[S(i)] + \sum_{j<i}c_j[S(j)] = [\Delta^{\sigma_i}(i+1)] + H_i[\Delta^{\sigma_i}(i)]\]
        for some integers $c_j \geq 0$.
        These equalities implies that there is an exact sequence
        \[0 \to \Delta^{\sigma_i}(i)^{H_i} \to \Delta(i+1) \to \Delta^{\sigma_i}(i+1) \to 0,\]
        and hence $(\mathcal{E}_i)$ holds as $\Delta^{\sigma_i}(i)=\Delta(i)$.
        On the other hand, assume $\overline{E_i}=0$.
        Then it suffices to show that the condition $(\overline{\mathcal{E}_i})$ holds if $D(A_A)$ has a $\nabla$-filtration by \Cref{AHLU}.
        Obviously, this is dually followed from the above discussion.

        Conversely, we prove the ``if'' part.
        Assume that $(\mathcal{E}_i)$ holds.
        We have $\Delta^{\sigma_i}(k)=\Delta(k)$ for $k \neq i+1$ again since $E_i=0$.
        By our assumption, we have an injection $\Delta(i)^{H_i} \to \Delta(i+1)$.
        Its image is equal to $\tr^{\sigma_i}_{i}\Delta(i+1)$ and the cokernel is isomorphic to 
        \[\begin{split}            
            \frac{\Delta(i+1)}{\tr^{\sigma_i}_{i}\Delta(i+1)} 
            &\cong \frac{\frac{P(i+1)}{\tr^e_{i+2}(P(i+1))}}{\tr^{\sigma_i}_{i}\left(\frac{P(i+1)}{\tr^e_{i+2}(P(i+1))}\right)}
            \cong \frac{\frac{P(i+1)}{\tr^e_{i+2}(P(i+1))}}{\frac{\tr^{\sigma_i}_{i}(P(i+1))}{\tr^e_{i+2}(P(i+1)) \cap \tr^{\sigma_i}_{i}(P(i+1))}}
            \cong \frac{\frac{P(i+1)}{\tr^e_{i+2}(P(i+1))}}{\frac{\tr^{\sigma_i}_{i}(P(i+1))+\tr^e_{i+2}(P(i+1))}{\tr^e_{i+2}(P(i+1))}} \\
            &\cong \frac{P(i+1)}{\tr^{\sigma_i}_{i}(P(i+1))+\tr^e_{i+2}(P(i+1))}
            = \frac{P(i+1)}{\tr^{\sigma_i}_{i}P(i+1)} \cong \Delta^{\sigma_i}(i+1),        
        \end{split}\]
        where the first and last isomorphisms are followed from \Cref{qh rmk} and the second one is from \Cref{qh rmk}.
        Thus $\Delta(i+1)$, and hence $_AA$, has a $\Delta^{\sigma_i}$-filtration.
        Dually, the condition $(\overline{\mathcal{E}_i})$ implies that $D(A_A)$ has a $\nabla$-filtration, and so $_AA$ has a $\Delta$-filtration by \Cref{AHLU} again.
        This completes the proof.
    \end{proof}
\end{thm}

The next corollary shows several properties of $i${\rm th}-twistable quasi-hereditary algebras.
\begin{cor}\label{ith-twist}
Let a quasi-hereditary algebra $(A,e)$ be $i${\rm th}-twistable.
Then the followings hold.
\begin{enumerate}[$(1)$]
\item $\dim\Hom_A(\Delta^{\sigma_i}(i+1),\Delta^{\sigma_i}(i))=E_i$.
\item $\dim\Hom_A(\nabla^{\sigma_i}(i),\nabla^{\sigma_i}(i+1))=\overline{E_i}$.
\item If $(\overline{\mathcal{E}_i})$ (resp. $(\mathcal{E}_i)$) in \Cref{ith-twistable} holds for $(A,e)$, then $(\mathcal{E}_i)$ (resp. $(\overline{\mathcal{E}_i})$) does for $(A,\sigma_i)$.
\end{enumerate}
\begin{proof}
    $(1)$ Assume that $E_i=0$. 
    Then $\Delta^{\sigma_i}(i)=\Delta(i)$, and hence $[\Delta^{\sigma_i}(i):S(i+1)]=0$.
    This shows that $\dim\Hom_A(\Delta^{\sigma_i}(i+1),\Delta^{\sigma_i}(i))=0$.
    Next, assume $E_i \neq 0$.
    Then there is an exact sequence
    \[0 \to \Delta(i+1)^{E_i} \to \Delta^{\sigma_i}(i) \to \Delta(i) \to 0,\]
    since $\Delta^{\sigma_i}(i) \cong P(i)/\tr^{\sigma_i}_{i+2}P(i) = P(i)/\tr^e_{i+2}P(i)$.
    Moreover our assumption implies $\overline{E_i}=0$, and hence $H_i=0$ by using \Cref{MO}.
    This gives the equality $\Delta^{\sigma_i}(i+1)=\Delta(i+1)$.
    By the exact sequence above, we get $\dim\Hom_A(\Delta^{\sigma_i}(i+1),\Delta^{\sigma_i}(i))=E_i$.

    $(2)$ This is similarly shown by arguments in the proof of $(1)$.

    $(3)$ By \Cref{ith-twistable}, the condition $(\mathcal{E}_i)$ or $(\overline{\mathcal{E}_i})$ holds for $(A,e)$. 
    From the above discussion, if $(A,e)$ satisfies $(\overline{\mathcal{E}_i})$, then $\Delta^{\sigma_i}(i+1)=\Delta(i+1)$ and there is an exact sequence
    \[0 \to \Delta(i+1)^{E_i} \to \Delta^{\sigma_i}(i) \to \Delta(i) \to 0.\]
    Since $\overline{E_i}=0$, we have 
    \[\dim\Ext^1_A(\Delta^{\sigma_i}(i+1),\Delta^{\sigma_i}(i))=\dim\Hom_A(\nabla^{\sigma_i}(i),\nabla^{\sigma_i}(i+1))=\overline{E_i}=0.\]
    Here the first equality uses \Cref{MO}.
    Moreover the above exact sequence shows that $\Delta^{\sigma_i}(i)$ has a submodule isomorphic to $\Delta^{\sigma_i}(i+1)^{E_i} = \Delta(i+1)^{E_i}$.
    Hence $(A,\sigma_i)$ satisfies $(\mathcal{E}_i)$.
    Dually if $(\mathcal{E}_i)$ holds for $(A,e)$, then $(\overline{\mathcal{E}_i})$ does so for $(A,\sigma_i)$.
\end{proof}
\end{cor}

\begin{eg}\label{eg twistable}
Consider a quiver $\xymatrix{1\ar[r]^\alpha \ar@/_10pt/[rr]_\beta & 2\ar[r]^\gamma &3 \ar[r]^\delta & 4}$ and an ideal $I=\langle \alpha\gamma\delta-\beta\delta \rangle$ of $KQ$, and put $A=KQ/I$. 
Now we have $24$ permutations on $\Lambda=\{1,2,3,4\}$.
In the following, we will write $\sigma_{(i_l \cdots i_2 i_1)}=\sigma_{(i_l, \dots, i_2, i_1)}$ as the product $\sigma_{i_l} \cdots \sigma_{i_2} \sigma_{i_1}$, where $i_k \in \{1,2,3\}$ for $1 \leq k \leq l$.
Let $\lambda=(i_l, \dots, i_2, i_1)$ be a sequence of elements of $\{1,2,3\}$ and $\Delta^{(\lambda)}$ be standard modules with respect to the $\sigma_{(\lambda)}$-order.
First, for each $(A,\sigma_{(\lambda)})$, we define the {\bfseries neighbor Hom-Ext${^1}$ biquiver} of $\Delta^{(\lambda)}$.
Its vertices $i$ are corresponding to standard modules $\Delta^{(\lambda)}(i)$.
For arrows, there are two kinds; solid and dotted.
If $\sigma_{(\lambda)}(j)-\sigma_{(\lambda)}(i) \neq 1$, there are no arrows from $i$ to $j$.
For $i,j$ satisfying $\sigma_{(\lambda)}(j)-\sigma_{(\lambda)}(i)=1$, solid arrows $i \to j$ are given if $\Ext^1_A(\Delta^{(\lambda)}(i),\Delta^{(\lambda)}(j)) \neq 0$, and dotted ones $i \dra j$ exist if $\Hom_A(\Delta^{(\lambda)}(i),\Delta^{(\lambda)}(j)) \neq 0$.
Similarly, the neighbor Hom-Ext${^1}$ biquiver of $\nabla^{(\lambda)}$ is also defined.
If $(A,\sigma_{(\lambda)})$ is quasi-hereditary, then the existence of arrows in neighbor Hom-Ext${^1}$ biquiver of $\nabla^{(\lambda)}$ are also determined by ones in the biquiver of $\Delta^{(\lambda)}$, using \Cref{MO,ith-twist}.
The following is the table of neighbor Hom-Ext${^1}$ biquivers of $\Delta^{(\lambda)}$ over the algebra $A$ for all possible $\lambda$.

\begin{multicols}{2}
    \begin{description}
    \item[$e$] $\xymatrix{1\ar[r]  & 2\ar[r] & 3\ar[r] & 4}$
    \item[$\sigma_{(1)}$] $\xymatrix{2\ar@{.>}[r] & {\begin{smatrix}1\\2\end{smatrix}}\ar[r] & 3\ar[r] & 4}$
    \item[$\sigma_{(2)}$] $\xymatrix{1\ar[r] & 3\ar@{.>}[r] & {\begin{smatrix}2\\3\end{smatrix}}\ar[r] & 4}$
    \item[$\sigma_{(3)}$] $\xymatrix{1\ar[r] & 2 & 4\ar@{.>}[r] & {\begin{smatrix}3\\4\end{smatrix}}}$
    \item[$\sigma_{(21)}$] $\xymatrix{2\ar[r] & 3\ar@{.>}[r] & {\begin{smatrix}&1\\2&&3\\3\end{smatrix}} \ar[r] & 4}$
    \item[$\sigma_{(31)}$] $\xymatrix{2\ar@{.>}[r] & {\begin{smatrix}1\\2\end{smatrix}} & 4\ar@{.>}[r] & {\begin{smatrix}3\\4\end{smatrix}}}$
    \item[$\sigma_{(12)}$] $\xymatrix{3\ar@{.>}[r] & {\begin{smatrix}1\\3\end{smatrix}}\ar[r] & {\begin{smatrix}2\\3\end{smatrix}}\ar[r] & 4}$
    \item[$\sigma_{(32)}$] $\xymatrix{1\ar[r] & 3\ar[r] & 4\ar@{.>}[r] & {\begin{smatrix}2\\3\\4\end{smatrix}}}$
    \item[$\sigma_{(23)}$] $\xymatrix{1 & 4 & 2\ar[r] & {\begin{smatrix}3\\4\end{smatrix}}}$
    \item[$\sigma_{(121)}$] $\xymatrix{3\ar@{.>}[r] & {\begin{smatrix}2\\3\end{smatrix}}\ar@{.>}[r] & {\begin{smatrix}&1\\2&&3\\3\end{smatrix}}\ar[r] & 4}$
    \item[$\sigma_{(3 2 1)}$] $\xymatrix{2\ar[r] & 3\ar[r] & 4\ar@{.>}[r] & {\begin{smatrix}&1\\2&&3\\3\\&4\end{smatrix}}}$
    \item[$\sigma_{(231)}$] $\xymatrix{2 & 4 & {\begin{smatrix}1\\2\end{smatrix}} & {\begin{smatrix}3\\4\end{smatrix}}}$
    \item[$\sigma_{(3 1 2)}$] $\xymatrix{3\ar@{.>}[r] & {\begin{smatrix}1\\3\end{smatrix}} & 4\ar@{.>}[r] & {\begin{smatrix}2\\3\\4\end{smatrix}}}$
    \item[$\sigma_{(232)}$] $\xymatrix{1 & 4\ar@{.>}[r] & {\begin{smatrix}3\\4\end{smatrix}}\ar@{.>}[r] & {\begin{smatrix}2\\3\\4\end{smatrix}}}$
    \item[$\sigma_{(123)}$] $\xymatrix{4 & 1\ar[r] & 2\ar[r] & {\begin{smatrix}3\\4\end{smatrix}}}$
    \item[$\sigma_{(3121)}$] $\xymatrix{3\ar@{.>}[r] & {\begin{smatrix}2\\3\end{smatrix}}\ar[r] & 4\ar@{.>}[r] & {\begin{smatrix}&1\\2&&3\\3\\&4\end{smatrix}}}$
    \item[$\sigma_{(2 3 2 1)}$] $\xymatrix{2 & 4\ar@{.>}[r] & {\begin{smatrix}3\\4\end{smatrix}}\ar@{.>}[r] & {\begin{smatrix}&1\\2&&3\\3\\&4\end{smatrix}}}$
    \item[$\sigma_{(1231)}$] $\xymatrix{4 & 2\ar@{.>}[r] & {\begin{smatrix}1\\2\end{smatrix}} & {\begin{smatrix}3\\4\end{smatrix}}}$
    \item[$\sigma_{(2312)}$] $\xymatrix{3\ar[r] & 4 & {\begin{smatrix}1\\3\end{smatrix}}\ar[r] & {\begin{smatrix}2\\3\\4\end{smatrix}}}$
    \item[$\sigma_{(1232)}$] $\xymatrix{4 & 1 & {\begin{smatrix}3\\4\end{smatrix}}\ar@{.>}[r] & {\begin{smatrix}2\\3\\4\end{smatrix}}}$
    \item[$\sigma_{(23121)}$] $\xymatrix{3\ar[r] & 4\ar@{.>}[r] & {\begin{smatrix}2\\3\\4\end{smatrix}}\ar@{.>}[r] & {\begin{smatrix}&1\\2&&3\\3\\&4\end{smatrix}}}$
    \item[$\sigma_{(12321)}$] $\xymatrix{4 & 2\ar[r] & {\begin{smatrix}3\\4\end{smatrix}}\ar@{.>}[r] & {\begin{smatrix}&1\\2&&3\\3\\&4\end{smatrix}}}$
    \item[$\sigma_{(12312)}$] $\xymatrix{4\ar@{.>}[r] & {\begin{smatrix}3\\4\end{smatrix}}\ar@{.>}[r] & {\begin{smatrix}1\\3\end{smatrix}}\ar[r] & {\begin{smatrix}2\\3\\4\end{smatrix}}}$
    \item[$\sigma_{(123121)}$] $\xymatrix{4\ar@{.>}[r] & {\begin{smatrix}3\\4\end{smatrix}}\ar@{.>}[r] & {\begin{smatrix}2\\3\\4\end{smatrix}}\ar@{.>}[r] & {\begin{smatrix}&1\\2&&3\\3\\&4\end{smatrix}}}$
    \end{description}
\end{multicols}

Clearly $(A,\sigma_{(123121)})$ is quasi-hereditary since all standard modules are projective and satisfy $[P(i):S(i)]=1$ for all $i \in \Lambda$.
By using $(\mathcal{E}_i)$ in \Cref{ith-twistable}, we have the following diagram which shows that if the source of an arrow gives a quasi-hereditary structure, then so does the target.
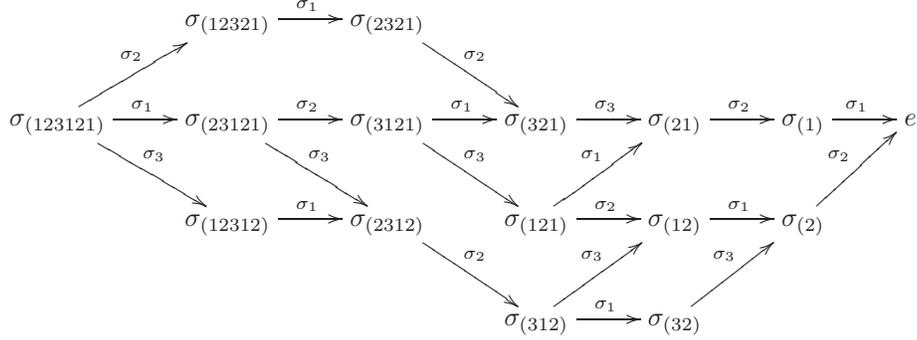
\[\xymatrix{
    & \sigma_{(12321)}\ar[r]^{\sigma_1} & \sigma_{(2321)}\ar[dr]^{\sigma_2} \\
    \sigma_{(123121)}\ar[r]^{\sigma_1}\ar[ur]^{\sigma_2}\ar[dr]^{\sigma_3} & \sigma_{(23121)}\ar[r]^{\sigma_2}\ar[dr]^{\sigma_3} & \sigma_{(3121)}\ar[r]^{\sigma_1}\ar[dr]^{\sigma_3} & \sigma_{(321)}\ar[r]^{\sigma_3} & \sigma_{(21)}\ar[r]^{\sigma_2} & \sigma_{(1)}\ar[r]^{\sigma_1} & e\\
    & \sigma_{(12312)}\ar[r]^{\sigma_1} & \sigma_{(2312)}\ar[dr]^{\sigma_2} & \sigma_{(121)}\ar[r]^{\sigma_2}\ar[ur]^{\sigma_1} & \sigma_{(12)}\ar[r]^{\sigma_1} & \sigma_{(2)}\ar[ur]^{\sigma_2}\\
    &&& \sigma_{(312)}\ar[r]^{\sigma_1}\ar[ur]^{\sigma_3} & \sigma_{(32)}\ar[ur]^{\sigma_3}
}\]

However applying \Cref{ith-twistable} to the quasi-hereditary algebra $(A,\sigma_{(12321)})$, we recognize that it is not $3${\rm rd}-twistable, i.e., $\sigma_{(1231)}$ does not give a quasi-hereditary structure of $A$.
Similarly, $(A,\sigma_{(12312)})$ is not $2${\rm nd}-twistable, and hence $\sigma_{(1232)}$ does not give a quasi-hereditary structure of $A$.
Focus on the two permutations $\sigma_{(1231)}$ and $\sigma_{(1232)}$.
Then we show the others which do not give quasi-hereditary structures. 
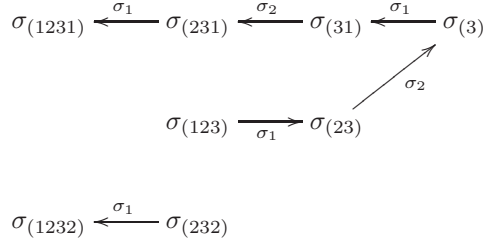
\[\xymatrix{
\sigma_{(1231)} & \sigma_{(231)}\ar[l]_{\sigma_1} & \sigma_{(31)}\ar[l]_{\sigma_2} & \sigma_{(3)}\ar[l]_{\sigma_1}\\
& \sigma_{(123)}\ar[r]_{\sigma_1} & \sigma_{(23)}\ar[ur]_{\sigma_2}\\
\sigma_{(1232)} & \sigma_{(232)}\ar[l]_{\sigma_1}
}\]
If some permutations in this diagram give quasi-hereditary structures, then $\sigma_{(1231)}$ or $\sigma_{(1232)}$ also does, a contradiction.
Hence the permutations in this diagram do not give quasi-hereditary structures of $A$.
Now we complete checking whether each permutation gives a quasi-hereditary structure or not for the algebra $A$.
\end{eg}

\section{Connectedness}\label{Connectedness}
In this section, we will argue the ``connectivity'' of quasi-hereditary structures.
In \Cref{eg twistable}, we got all permutations giving quasi-hereditary structures from one.
In general, we can obtain all those from one by checking repeatedly whether each quasi-hereditary algebra satisfies the condition $(\mathcal{E}_i)$ or $(\overline{\mathcal{E}_i})$.
Recall that we put $H_i=\dim\Hom_A(\Delta(i),\Delta(i+1))$ and $E_i=\dim\Ext_A^1(\Delta(i),\Delta(i+1))$.

\begin{df}\label{def of connected}
    Two permutations $\sigma$ and $\tau$ giving quasi-hereditary structures of $A$ are said to be {\bfseries connected} if the following condition holds:
    There is a decomposition $\tau \sigma^{-1} = \sigma_{i_l} \cdots \sigma_{i_1}$ into the product of adjacent transpositions such that all $\sigma_{i_k} \cdots \sigma_{i_1} \sigma$ for $1 \leq k \leq l$ also give quasi-hereditary structures.
\end{df}

If any two permutations giving quasi-hereditary structures are connected, we also say that the quasi-hereditary structures are connected. 
To show the connectedness of quasi-hereditary structures, we first prove that $e$ and another are connected, and in \Cref{connected} we prove that any two permutations are connected.


\begin{lem}\label{twist}
    Let $e,\sigma$ give quasi-hereditary structures of $A$ with $e \neq \sigma$.
    Then for the minimum element $i \in \Lambda$ satisfying $\sigma(i) > \sigma(i+1)$, it holds that $E_i H_i =0$.
    \begin{proof}
        Take the minimum element $i > 1$ with $\sigma(i) > \sigma(i+1)$.
        Then we have $\sigma(j) < \sigma(i)$ for $j < i$.
        Indeed, if there exists $j<i$ with $\sigma(j) > \sigma(i)$, some $j \leq k < i$ satisfies $\sigma(k+1) < \sigma(k)$.
        However this contradicts the minimality of $i$.
        So there is an epimorphism $\Delta^{\sigma}(i) \to \Delta(i)$.
        Additionally its kernel surjectively maps to $\Delta(i+1)^{E_i}$ since we have a commutative diagram 
        \[\xymatrix{
            0 \ar[r] & \kernel{f} \ar[r]\ar[d] & \Delta^{\sigma}(i) \ar[r]^-f\ar[d] & \Delta(i) \ar[r]\ar@{=}[d] & 0\\
            0 \ar[r] & \Delta(i+1)^{E_i} \ar[r] & P(i)/\tr^{e}_{i+2}P(i) \ar[r] & \Delta(i) \ar[r] & 0,
        }\]
        where each row is exact and the middle vertical homomorphism is surjective.
        So we have\\ $[\Delta^{\sigma}(i):S(i)] \geq 1 + E_i H_i$.
        Since $(A,\sigma)$ is quasi-hereditary, the product $E_i H_i$ must be zero.
    \end{proof} 
\end{lem}

\begin{lem}\label{mini element}
    Let $e,\sigma$ give quasi-hereditary structures of $A$ with $e \neq \sigma$.
    Then for the minimum element $i \in \Lambda$ satisfying $\sigma(i) > \sigma(i+1)$, it gives quasi-hereditary structure.
    \begin{proof}
        By \Cref{MO,twist}, we have $E_i=0$ or $\overline{E_i}=0$.
        It is enough to show that if $E_i=0$, then the condition $(\mathcal{E}_i)$ holds for $(A,e)$.
        Indeed, we can dually show that $(\overline{\mathcal{E}_i})$ holds if $\overline{E_i}=0$. 
        Define $i_{\sigma+}$ by $\sigma(i_{\sigma+})=\sigma(i)+1$.
        We have already shown that $\sigma(j) < \sigma(i)$ for $j < i$ in the proof of \Cref{twist}.

        So there are canonical epimorphisms $p:P(i+1)/\tr^{\sigma}_{i_{\sigma+}}(P(i+1)) \to \Delta(i+1)$ and $q:\Delta^\sigma(i) \to \Delta(i)$ with kernels $\tr^e_{i+2}\left(P(i+1)/\tr^{\sigma}_{i_{\sigma+}}(P(i+1))\right)$ and $\tr^e_{i+1}(\Delta^\sigma(i)^{H_i})$, respectively.
        Moreover the direct summands of the top of $\tr^e_{i+1}(\Delta^\sigma(i)^{H_i})$ do not include $S(i+1)$ since $E_i=0$.        
        Hence we have $\tr^e_{i+1}(\Delta^\sigma(i)^{H_i})=\tr^e_{i+2}(\Delta^\sigma(i)^{H_i})$.

        Consider the set
        \[\left\{f_k:\Delta^\sigma(i) \overset{\iota_k}{\to} \Delta^\sigma(i)^{[P(i+1):\Delta^\sigma(i)]} \to \frac{P(i+1)}{\tr^{\sigma}_{i_{\sigma+}}(P(i+1))} \middle| \iota_k:k\text{-th} \text{ inclusion}\right\}_{1 \leq k \leq [P(i+1):\Delta^\sigma(i)]}.\]
        Then this forms a basis of $\Hom_A(\Delta^\sigma(i),P(i+1)/\tr^{\sigma}_{i_{\sigma+}}(P(i+1)))$, since $P(i+1)/\tr^{\sigma}_{i_{\sigma+}}(P(i+1))$ has a submodule isomorphic to $\Delta^\sigma(i)^{[P(i+1):\Delta^\sigma(i)]}$ and $[P(i+1)/\tr^{\sigma}_{i_{\sigma+}}(P(i+1)):S(i)] = [P(i+1):\Delta^\sigma(i)]$.
        On the other hand, the induced map 
        \[\Hom_A(\Delta^\sigma(i),p) : \Hom_A\left(\Delta^\sigma(i),\frac{P(i+1)}{\tr^{\sigma}_{i_{\sigma+}}(P(i+1))}\right) \to \Hom_A(\Delta^\sigma(i),\Delta(i+1)); f \mapsto p f\]
        is surjective.
        Indeed, we have the exact sequence
        {\small
        \[\begin{split}
            0 &\to \Hom_A\left(\Delta^\sigma(i),\tr^e_{i+2}\left(\frac{P(i+1)}{\tr^{\sigma}_{i_{\sigma+}}(P(i+1))}\right)\right)\\ &\to \Hom_A\left(\Delta^\sigma(i),\frac{P(i+1)}{\tr^{\sigma}_{i_{\sigma+}}(P(i+1))}\right) \to \Hom_A(\Delta^\sigma(i),\Delta(i+1))
        \end{split}\]
        }
        whose dimension of each non-zero term is $\left[\tr^e_{i+2}\left(\frac{P(i+1)}{\tr^{\sigma}_{i_{\sigma+}}(P(i+1))}\right):S(i)\right]$, $[P(i+1):\Delta^\sigma(i)]$, and $H_i$, respectively.
        Further we can easily check that 
        \[\begin{split}
            [P(i+1):\Delta^\sigma(i)]
            &=\left[\frac{P(i+1)}{\tr^{\sigma}_{i_{\sigma+}}(P(i+1))}:S(i)\right]\\
            &=\left[\tr^e_{i+2}\left(\frac{P(i+1)}{\tr^{\sigma}_{i_{\sigma+}}(P(i+1))}\right):S(i)\right]+[\Delta(i+1):S(i)]\\
            &=\left[\tr^e_{i+2}\left(\frac{P(i+1)}{\tr^{\sigma}_{i_{\sigma+}}(P(i+1))}\right):S(i)\right]+H_i,
        \end{split}\]
        and this shows that $\Hom_A(\Delta^\sigma(i),p)$ is surjective.
        Thus, changing suffixes if necessary we let $\{pf_1,\dots,pf_{H_i}\}$ be a basis of $\Hom_A(\Delta^{\sigma}(i),\Delta(i+1))$.
        Now put 
        \[f=(f_1,\dots,f_{H_i}):\Delta^{\sigma}(i)^{H_i} \to \frac{P(i+1)}{\tr^{\sigma}_{i_{\sigma+}}(P(i+1))},\]
        then it is an injection.

        Then by \Cref{trace}, for some $g : \Delta(i)^{H_i} \to \Delta(i+1)$, we have the following commutative diagram
        \[\xymatrix{
            0 \ar[r] & \tr^e_{i+2}(\Delta^\sigma(i)^{H_i}) \ar[r]\ar[d]^-{\tr^e_{i+2}(f)} & \Delta^\sigma(i)^{H_i} \ar[r]^-q \ar[d]^-f & \Delta(i)^{H_i} \ar[r]\ar[d]^-g & 0\\
            0 \ar[r] & \tr^e_{i+2}\left(\frac{P(i+1)}{\tr^{\sigma}_{i_{\sigma+}}(P(i+1))}\right) \ar[r] & \frac{P(i+1)}{\tr^{\sigma}_{i_{\sigma+}}(P(i+1))} \ar[r]^-p & \Delta(i+1) \ar[r] & 0
        }\]
        where each row is exact.
        Applying the snake lemma to this diagram, we obtain the exact sequence 
        \[\kernel f \to \kernel g \to \Coker(\tr^e_{i+2}(f)) \to \Coker f.\]
        Immediately we get $\kernel f=0$ since $f$ is injective.
        Next, we will prove that there is an isomorphism $\Coker(\tr^e_{i+2}(f)) \cong \tr^e_{i+2}(\Coker(f))$.
        To get this isomorphism, we show the equality
        \[\image \tr^e_{i+2}(f) = \image f \cap \tr^e_{i+2}\left(\frac{P(i+1)}{\tr^{\sigma}_{i_{\sigma+}}(P(i+1))}\right).\]
        First, we have that $\image \tr^e_{i+2}(f) \subset \image f \cap \tr^e_{i+2}\left(\frac{P(i+1)}{\tr^{\sigma}_{i_{\sigma+}}(P(i+1))}\right)$, in general.
        Since we have equalities $\image\tr^e_{i+2}(f) = \tr^e_{i+2}(\Delta^\sigma(i)^{H_i}) = \tr^e_{i+1}(\Delta^\sigma(i)^{H_i})$, hence for $k \geq i+1$, we have
        \[0 = \left[\frac{\image{f}}{\image(\tr^e_{i+2}{f})}:S(k)\right] \geq \left[\frac{\image{f}\cap\tr^e_{i+2}\left(\frac{P(i+1)}{\tr^{\sigma}_{i_{\sigma+}}(P(i+1))}\right)}{\image(\tr^e_{i+2}{f})}:S(k)\right].\]
        On the other hand, by the definition of $f$,
        \[\left[\image{f}\cap\tr^e_{i+2}\left(\frac{P(i+1)}{\tr^{\sigma}_{i_{\sigma+}}(P(i+1))}\right):S(i)\right] = 0.\]
        Since the compositions of the top of $\image{f}\cap\tr^e_{i+2}\left(\frac{P(i+1)}{\tr^{\sigma}_{i_{\sigma+}}(P(i+1))}\right)$ are direct sum of $S(j)$'s for $j \geq i$, we conclude that 
        \[\image \tr^e_{i+2}(f) = \image f \cap \tr^e_{i+2}\left(\frac{P(i+1)}{\tr^{\sigma}_{i_{\sigma+}}(P(i+1))}\right).\]
        So, by \Cref{trace} there are isomorphisms
        \[
            \Coker(\tr^e_{i+2}(f)) 
            = \frac{\tr^e_{i+2}\left(\frac{P(i+1)}{\tr^{\sigma}_{i_{\sigma+}}(P(i+1))}\right)}{\image f \cap \tr^e_{i+2}\left(\frac{P(i+1)}{\tr^{\sigma}_{i_{\sigma+}}(P(i+1))}\right)}
            \cong \tr^e_{i+2}\left(\frac{\frac{P(i+1)}{\tr^{\sigma}_{i_{\sigma+}}(P(i+1))}}{\image f}\right)
            = \tr^e_{i+2}(\Coker(f)).
        \]
        Moreover we will show that the given map $\Coker(\tr^e_{i+2}(f)) \to \Coker f$ is passed through\\ $\tr^e_{i+2}(\Coker(f))$ with an isomorphism $\Coker(\tr^e_{i+2}(f)) \cong \tr^e_{i+2}(\Coker(f))$.
        For simplicity, we put $\Delta^\sigma(i)^{H_i} = X$ and $\frac{P(i+1)}{\tr^{\sigma}_{i_{\sigma+}}(P(i+1))} = Y$.
        Now we consider a diagram
        \[\xymatrix{
            \tr^e_{i+2}X \ar[r]^-{\tr^e_{i+2}(f)}\ar@{=}[d] & \tr^e_{i+2}Y \ar[r]^-{\pi'}\ar@{=}[d] & \Coker(\tr^e_{i+2}{f})\ar[d]^-{\varphi}\\
            \tr^e_{i+2}X \ar[r]^-{\tr^e_{i+2}(f)}\ar[d] & \tr^e_{i+2}Y \ar[r]^-{\tr^e_{i+2}(\pi)}\ar[d] & \tr^e_{i+2}(\Coker{f})\ar[d]\\
            X \ar[r]^-f & Y \ar[r]^-{\pi} & \Coker{f}.
        }\]
        Then upper and lower diagrams are commutative by the universality of cokernels and \Cref{trace}, respectively.
        Moreover since $\pi$ is surjective, so is $\varphi \pi' = \tr^e_{i+2}(\pi)$ by \Cref{trace} again.
        So $\varphi$ is an epimorphism, and hence it must be isomorphic as the argumentation above.
        Thus the composition map $\Coker(\tr^e_{i+2}(f)) \overset{\varphi}{\to} \tr^e_{i+2}(\Coker(f)) \hookrightarrow \Coker f$ is injective.
        Summing up, we have $\kernel g=0$ i.e. $g:\Delta(i)^{H_i} \to \Delta(i+1)$ is a monomorphism.
        This shows that $\Delta(i)^{H_i}$ is isomorphic to a submodule of $\Delta(i+1)$, and hence the condition $(\mathcal{E}_i)$ holds.
    \end{proof}
\end{lem}

Moreover we similarly get the following.
\begin{cor}\label{mini element sigma}
    Let $e,\sigma$ give quasi-hereditary structures of $A$ with $e \neq \sigma$.
    Then for the minimum element $i \in \Lambda$ with respect to the $\sigma$-order satisfying $i > i_{\sigma+} = \sigma^{-1}(\sigma(i)+1)$, the permutation $\sigma \sigma_i$ gives quasi-hereditary structure.
\end{cor}

\begin{prop}\label{connected induction}
    Let $e,\sigma$ give quasi-hereditary structures of $A$ with $e \neq \sigma$.
    Then there is a minimal decomposition $\sigma = \sigma_{i_l} \cdots \sigma_{i_1}$ into the product of adjacent transpositions such that $\sigma \sigma_{i_l}$ gives a quasi-hereditary structure.
    \begin{proof}
        Take $i_l$ as the $i$ given in \Cref{mini element sigma}.
        Then $\sigma \sigma_{i_l}$ gives a quasi-hereditary structure.
    \end{proof}
\end{prop}

By using induction on the length of $\sigma$, we conclude the following.

\begin{thm}\label{e connected}
    Let $e,\sigma$ give quasi-hereditary structures.
    Then they are connected.
    \begin{proof}
        We show this statement by induction on the minimal length $l$ of the decomposition of $\sigma$ into the product of adjacent transpositions.
        For $l=1$, it is clear.
        Now we show the claim for $l>1$.
        By \Cref{connected induction}, we showed that $(A,\sigma_{i_l}\sigma)$ is quasi-hereditary for some $\sigma_{i_l}$ which is the leftmost in a minimal decomposition of $\sigma$.
        Now $\sigma_{i_l}\sigma$ can be written by a decomposition of length $l-1$, and hence we can use our induction hypothesis.
    \end{proof}
\end{thm}

Finally we show that any two permutations giving quasi-hereditary structures are connected by renaming the indices of idempotents and reducing to \Cref{e connected}.

\begin{thm}\label{connected}
    Any two permutations giving quasi-hereditary structures are connected.
    \begin{proof}
        Assume that $\sigma$ and $\tau$ give quasi-hereditary structures of an algebra $A$ with idempotents $e_1,\dots,e_n$.
        Then for the algebra $A'=A$ with idempotents $e'_1,\dots,e'_n$ where $e'_i=e_{\sigma^{-1}(i)}$ for each $i\in\Lambda$, permutations $e$ and $\tau\sigma^{-1}$ give quasi-hereditary structures of $A'$.
        By \Cref{e connected}, there exists a decomposition 
        \[\tau\sigma^{-1} = \sigma_{i_l} \cdots \sigma_{i_1}\]
        into the product of adjacent transpositions such that all $\sigma_{i_k} \cdots \sigma_{i_1}$ for $1 \leq k \leq l$ give quasi-hereditary structures of $A'$.
        This shows us that all $\sigma_{i_k} \cdots \sigma_{i_1} \sigma$ for $1 \leq k \leq l$ give quasi-hereditary structures of $A$, that is, $\sigma$ and $\tau$ are connected with respect to giving quasi-hereditary structures of $A$.
    \end{proof}
\end{thm}

Moreover, for two quasi-hereditary structures, we get a sequence of adjacent transpositions which induce the connectedness of them, by \Cref{connected induction}.
In particular, this sequence is determined by only the permutations and does not depend on the algebra.
\begin{cor}\label{decomp of connected}
    Let $\sigma,\tau$ give quasi-hereditary structures with $\sigma \neq \tau$.
    For $k = 1,2,\dots$, inductively take a minimal element $i_{k}$ with respect to the $(\sigma_{i_{k-1}} \cdots \sigma_{i_1} \sigma)$-order satisfying $\sigma_{i_{k-1}} \cdots \sigma_{i_1} \sigma(i_{k}) \neq n$ and
    \[\tau(i_{k}) > \tau\sigma^{-1} \sigma_{i_1} \cdots \sigma_{i_{k-1}}(\sigma_{i_{k-1}} \cdots \sigma_{i_1} \sigma (i_{k})+1).\]
    We take $i_1,i_2,\dots,i_k$ until those elements satisfying the above exist.
    If there is no $i_{k+1}$ satisfying the above, then we do not take $i_{k+1}$ and put $l=k$.
    Then the product $\sigma_{i_l} \cdots \sigma_{i_1}$ is a decomposition of $\tau\sigma^{-1}$ inducing the connectedness of $\sigma$ and $\tau$.
\end{cor}

\end{document}